\documentclass[12pt,twoside,final,psamsfonts]{amsart}

\usepackage[psamsfonts]{amssymb}
\usepackage{times,a4wide}

\usepackage{color}		% Need the color package
\usepackage{epsfig}

\theoremstyle{plain}
\newtheorem*{theorem*}{Theorem}

\newtheorem{theorem}{Theorem}[section]
\newtheorem{proposition}[theorem]{Proposition}
\newtheorem*{proposition*}{Proposition}
\newtheorem{corollary}[theorem]{Corollary}
\newtheorem*{corollary*}{Corollary}

\newtheorem*{lemma*}{Lemma}

\theoremstyle{definition}
\newtheorem{remark}[theorem]{Remark}
\newtheorem*{remark*}{Remark}
\newtheorem*{questions*}{Question}

\theoremstyle{definition}

\newtheorem*{definition*}{Definition}
\newcommand{\nc}{\newcommand}

\newcommand{\N}{{\mathbb N}}
\newcommand{\DD}{{\mathbb D}}
\newcommand{\R}{{\mathbb R}}

\newcommand{\T}{{\mathbb{T}}}

\newcommand{\Int}{\operatorname{Int}}
\newcommand{\Hol}{\operatorname{Hol}}
\newcommand{\Har}{\operatorname{Har}}

\nc{\Intl}{\Int_{l^{\infty}}}

\newcommand{\lra}{\longrightarrow}
\newcommand{\lmto}{\longmapsto}
\newcommand{\eps}{\varepsilon}
\newcommand{\vp}{\varphi}
\nc{\bea}{\begin{eqnarray}}
\nc{\eea}{\end{eqnarray}}
\nc{\beqa}{\begin{eqnarray*}}
\nc{\eeqa}{\end{eqnarray*}}
\nc{\Hi}{H^{\infty}}
\nc{\loi}{\ell^{\infty}}
\nc{\NL}{N^+\vert \Lambda}
\nc{\liL}{\lambda\in\Lambda}
\nc{\nn}{\nonumber}
\nc{\hf}{{\mathcal H}_{\Phi}}
\nc{\hF}{{\mathcal H}_{\Phi}}

\newenvironment{proof*}{\vskip 2mm\noindent {}}{$\blacksquare$ \vskip 2mm}
\numberwithin{equation}{section}

\newcommand{\eit}{e^{i\theta}}

\nc{\card}{\operatorname{card}}
\nc{\tsn}{\tilde{\sigma}_n}
\nc{\tsk}{\tilde{\sigma}_k}
\nc{\tskp}{\tilde{\sigma}_k^+}
\nc{\tskm}{\tilde{\sigma}_k^-}
\nc{\dst}{\displaystyle}
\nc{\vt}{\vartheta}

\renewcommand{\qedsymbol}{$\blacksquare$}

\title[Interpolation and harmonic majorants
in big Hardy-Orlicz spaces]{Interpolation and harmonic marjorants in 
big Hardy-Orlicz spaces}

\author{Andreas Hartmann}

\address{Laboratoire de Math\'ematiques Pures de Bordeaux,
Universit\'e Bordeaux I, 351 cours de la Lib\'eration,
33405 Talence, France}
\email{Andreas.Hartmann@math.u-bordeaux1.fr}

\date{\today}



\begin{document}

\begin{abstract}
Free interpolation in Hardy spaces is caracterized by the
well-known Carleson condition. The result extends to Hardy-Orlicz
spaces contained in the scale of classical Hardy spaces $H^p$, $p>0$. 
For the Smirnov and the  Nevanlinna classes, interpolating sequences 
have been characterized in a recent paper in terms of the existence
of harmonic majorants 
(quasi-bounded in the case of the Smirnov class). 
Since the Smirnov class can be regarded as the
union over all Hardy-Orlicz spaces associated with a so-called
strongly convex function, it is natural to ask 
how the condition changes from the Carleson condition in 
classical Hardy spaces to harmonic majorants in the Smirnov 
class. The aim of this paper is to narrow down this gap from
the Smirnov class to ``big'' Hardy-Orlicz spaces.
More precisely, we characterize interpolating
sequences for a class of Hardy-Orlicz spaces that
carry an algebraic structure and that are strictly bigger than $\bigcup_{p>0} 
H^p$. It turns out that the interpolating sequences are again
characterized by the existence of quasi-bounded majorants, but
now the weights of the majorants have to be in suitable
Orlicz spaces. The existence of harmonic majorants in such
Orlicz spaces will also be discussed in the general situation. 
We finish the paper with an example of a separated Blaschke
sequence that is interpolating for certain Hardy-Orlicz spaces
without being interpolating for slightly smaller ones.
\end{abstract}

\maketitle

\section{Introduction}

For a sequence $\Lambda$ in the unit disk $\DD$ and a space
of holomorphic functions $X\subset \Hol(\DD)$, the interpolation
problem consists in describing the trace of $X$ on $\Lambda$,
i.e.\ the set of restrictions $X|\Lambda$ regarded as a sequence
space. This problem has been considered for many spaces and also for
domains different from $\DD$.

In this paper, which can be regarded as a continuation of the
work in \cite{HMNT},
we will focus on spaces included in the Nevanlinna or the
Smirnov class, and in particular on Hardy-like spaces. Carleson
described in 1958 the interpolating sequences for the
Hardy space of bounded analytic functions on the unit disk by the
condition that is now called the Carleson condition \cite{carl}. 
It turns out that this condition still charactizes interpolating sequences
in a much broader situation.
It was successively proved to be the right condition  
in $H^p$, $p\ge 1$ \cite{ShHSh},
in $H^p$, $p < 1$ \cite{Kab}, and for Hardy-Orlicz spaces
$\hf$ contained in the scale $H^p$, $p>0$, in \cite{Har}
(even for the weaker notion of {\it free} interpolation, see
the definition below).

A natural question is then to ask what happens beyond $\bigcup_{p>0}H^p$?
It is here where we enter into the territory of ``big'' Hardy-Orlicz spaces.

Two observations should be made at this junction.
First, it is known on the one hand
that the union of all Hardy-Orlicz spaces corresponds to 
the Smirnov class (see \cite{RR}), and so especially the union of
big Hardy-Orlicz spaces. On the other hand, interpolating sequences
for the Smirnov class (and the Nevanlinna class) have been characterized
recently in the paper \cite{HMNT}. So, one could try to examine
interpolating sequences for big Hardy-Orlicz spaces in the light
of these results. Since the characterization for the Smirnov and
Nevanlinna classes is no longer given by the Carleson condition 
this leads to the still open question where exactly the Carleson condition
ceases to be valid (see also the Question at the end of this section).

The second observation concerns the notion of interpolating
sequences itself. When one wants to characterize for exemple 
interpolating sequences for the Smirnov or the Nevanlinna class,
which kind of natural trace space can one hope for a priori?
The problem had been studied in the past for a priori fixed
trace spaces (see \cite{Na56} for the Nevanlinna class and
\cite{yana2} for the Smirnov class), but those traces turned out
to be too big (Naftalevi\v c's case) or too small (Yanagihara's case),
see \cite{HMNT} for more detailed comments on this. The notion
that finally appeared to be natural is that of {\it free} interpolating
sequences. This notion goes back to work by Vinogradov and Havin
%PROBABLY RATHER VINOGRADOVS '76 paper (... in spaces with L^p norm)
%then Vinogradov-Havin 74 where the interpolation remains in H^infty
in the middle of the seventies
and has been used for general interpolation problems
by Nikolski and Vasyunin for Hilbert spaces (it works also
in certain Banach spaces). It has been succesfully used in \cite{HMNT} for the
Smirnov class and the Nevanlinna class (the latter being even not
a topological vector space).

%{\sc Check Vinogradov Havin!!!}

Let us introduce the notion of free interpolation that we will
use in this paper. %(going back to that of Vasyunin and Nikolski).

\begin{definition*}\label{ideal} A sequence space $l$ is called {\it ideal} if
$\ell^\infty l\subset l$, i.e.\ whenever $(a_n)_n\in l$ and  $ (\omega_n)_n\in
\ell^\infty$,  then also $(\omega_n a_n)_n\in l$.
\end{definition*}

\begin{definition*}\label{freeint} Let $X$ be a space of holomorphic functions
in $\DD$. A sequence $\Lambda\subset \DD$ is called {\it free interpolating for
$X$} if $X\vert \Lambda$ is ideal. We denote $\Lambda\in \Int X$.
\end{definition*}

\begin{remark}
\label{l-infinit}
For any function algebra $X$ containing the constants,
$\Lambda\in\Int X$ --- or equivalently
$X\vert\Lambda$ is ideal --- 
if and only if
\beqa
      \loi\subset X\vert\Lambda.
\eeqa

This remark is quite easy to check and a proof is given in \cite{HMNT}.
Let us repeat this proof here for completeness.
The inclusion is obviously necessary. In order to see that it is sufficient
notice that, by assumption, for any $(\omega_\lambda)_\lambda\in \loi$ there
exists $g\in X$ such that $g(\lambda)=\omega_\lambda$. Thus, if
$(f(\lambda))_\lambda\in X\vert \Lambda$, then the sequence of values
$(\omega_\lambda f(\lambda))_\lambda$ can be interpolated by  $fg\in X$.

An almost trivial but nevertheless useful remark in our
context is that if 
$\Lambda\in \Int X$ (i.e.\ $X\vert\Lambda$ is ideal) then $\Lambda'\in \Int X$
(i.e.\ $X\vert\Lambda'$ is ideal) for any subsequence
$\Lambda'\subset\Lambda$.
\end{remark}

%{\sc Definition of free interpolation : put the almost trivial remark
%here}

In the case of the Nevanlinna and Smirnov classes the characterization
of free interpolating sequences is given in terms of 
%will be expressed in terms of 
the existence of (quasi-bounded) harmonic majorants of a certain 
density %function 
associated with %the Blaschke product of 
the sequence $\Lambda$. This density is expressed in
terms of %the logarithm of
Blaschke products.
More precisely, let $b_{\lambda}(z)=\frac{|\lambda|}{\lambda}
\frac{\lambda-z}{1-\overline{\lambda}z}$ be the elementary Blaschke
factor (or M\"obius transform). For a sequence $\Lambda$ 
satisfying the Blaschke condition $\sum_{\lambda\in\Lambda} (1-|\lambda|^2)
<\infty$ --- which will be assumed throughout this paper --- we set
$B=B_{\Lambda}=\prod_{\lambda\in\Lambda}b_{\lambda}$ for the
corresponding Blaschke product, and $B_{\lambda}:=B_{\Lambda\setminus
\{\lambda\}}$. Define then
\beqa
 \varphi_\Lambda (z) :=
 \begin{cases}
 \log |B_\lambda(\lambda)|^{-1}\quad&\textrm{if $z=\lambda \in \Lambda$}\\
 \ 0 \quad&\textrm{if $z \notin \Lambda$}.
 \end{cases}
\eeqa

The result of \cite{HMNT} we are interested in is the following
(we give here only an almost complete form).
Recall that the Smirnov class $N^+$ is the class of holomorphic 
functions $f$ on the unit disk such that
\beqa
 \lim_{r \to 1}\frac{1}{2\pi}\int_{0}^{2\pi}
       \log^+|f(re^{i\theta})|\;d\theta=\frac{1}{2\pi}\int_{0}^{2\pi}
       \log^+|f(e^{i\theta})|\;d\theta.
\eeqa
(Note that the existence of the limit on the left hand side implies
the existence of the boundary values of $f$ on $\T$ a.e. appearing
on the right hand side).
%that we cite here in 
%(we will not state the complete theorem here).
 \begin{theorem*}[\cite{HMNT}]%\label{thmCNS}
Let $\Lambda$ be a sequence in $\DD$. The following statements are
equivalent:
\begin{itemize}
\item[(a)] $\Lambda$ is a free
interpolating sequence for the Smirnov class $N^+$: $\Lambda\in
\Int \hf$.

\item[(b)] $\varphi_\Lambda$ admits a quasi-bounded harmonic
majorant, i.e.\ there exists a positive weight $w\in L^1(\T)$
such that
\beqa
 \varphi_{\Lambda}(\lambda)\le P[w](\lambda)
 =\int_{\T} P_{\lambda}(\zeta) w(\zeta)\,dm(\zeta)
 =\int_{\T}\frac{1-|\lambda|^2}{|\zeta-{\lambda}|^2}
 w(\zeta)dm(\zeta),
\eeqa
where $P_{\lambda}(\zeta)=(1-|\lambda|^2)/|\zeta-\lambda|^2$ denotes 
the Poisson kernel.
\item[(c)] The trace space is given by
\[
N^+|\Lambda=l_{N^+}:=\{(a_{\lambda})_{\lambda} : \ \exists\
h\in\Har_+(\DD)\
\text{ quasi-bounded},\ 
h(\lambda)\ge \log^+|a_{\lambda}|, \ \liL \}.
\]
\item[(d)] $\lim_{n\to\infty}\sup_{(c_{\lambda)}\in {\mathcal
      B}_{\Lambda}} \sum_{\lambda:\vp_{\Lambda}(\lambda)\ge n}
  c_{\lambda} \vp_{\Lambda}(\lambda)=0$, where 
${\mathcal B}_{\Lambda}=\{(c_{\lambda}):c_{\lambda}\ge 0$ for any
$\lambda\in\Lambda$ and $\|\sum c_{\lambda}P_{\lambda}\|_{\infty}\le 1\}$.
\end{itemize}
\end{theorem*}

The condition (d) may appear quite technical. It is in a sense
a ``little-o''-version of the condition that characterizes
interpolation in the Nevanlinna class. In that case (see
\cite{HMNT}) free interpolating
sequences are characterized by the fact that $\vp_{\Lambda}$ admits
a positive harmonic majorant (be it quasi-bounded or not). Such
a majorant exists if and only if there is a constant $C$
such that for every finite positive sequence $(c_{\lambda})$ we have
\bea\label{charNevInt}
 \sum c_{\lambda} \vp_{\Lambda}(\lambda)\le C\|\sum c_{\lambda}
 P_{\lambda}\|_{\infty}.
\eea
(In view of discussions to come
we observe that $L^{\infty}=(L^1)^*$.)

In the light of the above theorem, and since the Hardy-Orlicz spaces we
are interested in are in a sense close to the Smirnov class,
% (they
%carry an algebraic structure and are much bigger than $H^p$-spaces),
it seems natural to seek for a condition in the spirit of condition
(b) in the theorem. The modification should involve a
more precise hypothesis on the weight $w$ adapted to the
defining function of the Hardy-Orlicz space that we will introduce
below. We will also see that the dual condition (d) of the theorem --- or
rather \eqref{charNevInt} --- has a counterpart in the Hardy-Orlicz
situation (replacing $(L^1)^*$ by the dual of a suitable Orlicz space).

%We now turn to Orlicz and Hardy-Orlicz spaces.

Let $\vp:{\R}\lra [0,\infty)$ be a convex, nondecreasing
function satisfying
\begin{itemize}
\item[(i)] $\lim_{t\to \infty} \vp (t)/t =\infty$
\item[(ii)] $\tilde{\Delta}_2$-condition: $\vp (t+2) \le M \vp (t) +K$,
$t \ge t_0$ for some constants $M,K \ge 0$ and $t_0\in \R$.
\end{itemize}
Such a function is called strongly convex (see \cite{RR}),
and one can associate with it the corresponding \emph{Hardy-Orlicz
class}
\beqa
      \mathcal{H}_{\vp\circ\log} 
 =\{f\in N^+:\int_\T\vp (\log |f(\zeta)|)\,d\sigma(\zeta)
      <\infty\},
\eeqa
where $f(\zeta)$ is the non-tangential boundary value of $f$ at
$\zeta\in {\T}$, which exists almost everywhere since $f\in N^+$.
Throughout this paper we shall use the notation
\beqa
 \Phi=\vp\circ\log.
\eeqa
It should be noted that the $\tilde{\Delta}_2$-condition is
formulated in such a way that $\Phi$ %=\vp\circ\log$ 
satisfies the
usual $\Delta_2$-condition: there exist constants $M',K'\ge 0$ and
$s_0$ such that for all $s\ge s_0$ we have
\bea\label{Delta}
 \Phi(2s)\le M'\Phi(s)+K'
\eea
(see also Section \ref{HO} for more on Hardy-Orlicz spaces).

In \cite{Har}, the following result was proved.

\begin{theorem*}[\cite{Har}]
Let $\vp$ be a strongly convex function satisfying (i), (ii) and
the $V_2$-condition:
\beqa
      2\vp(t)\le \vp(t+\alpha), \quad t\ge t_1,
\eeqa
where $\alpha>0$ is a suitable constant and $t_1\in\R$. Then $\Lambda
\subset {\DD}$ is free interpolating for $\hf$
if and only if
$\Lambda$  %is an interpolating sequence for $\Hi$, i.e.
satisfies the so-called Carleson condition:
\beqa
 \inf_{\lambda\in\Lambda}|B_{\lambda}(\lambda)|=\delta>0.
\eeqa
And in this case
\beqa
%    \mathcal{H}_{\vp\circ\log} 
 \hf\vert \Lambda
 =\{a=(a_\lambda)_\lambda: |a|_{\varphi}=
      \sum_{\liL} (1-|\lambda|)
      \varphi(\log |a_\lambda|)<\infty\}.
\eeqa
\end{theorem*}

Observe that the Carleson condition % of being interpolating for $\Hi$
can be reformulated in terms of $\vp_{\Lambda}$:
\bea\label{Intcst}
 M=\sup_{\lambda} \vp_{\Lambda}(\lambda) <+\infty,
\eea
i.e.\ $\vp_{\Lambda}$ is bounded and admits a fortiori a
harmonic majorant. We shall occasionally call $M$ the constant
associated with a sequence verifying the Carleson condition
(such sequences are usually called $\Hi$-interpolating sequences).

The conditions on $\vp$ in the theorem imply that %for all $\hf$
there exist $p,q\in (0,\infty)$ such that $H^p\subset \hf\subset
H^q$.
In particular, the $V_2$-condition
implies the inclusion $\hf\subset H^p$ for some $p>0$. This
$V_2$-condition has
a strong topological impact on the spaces. In fact, it guarantees that
metric
bounded sets are also bounded in the topology of the space (and so the
usual
functional analysis tools still apply in this situation; see
\cite{Har} for more on this and for further references).
It was not clear whether this was only a technical problem or if there
existed a critical growth for $\vp$ (below exponential growth
$\vp(t)=e^{pt}$ corresponding to $H^p$ spaces) giving a breakpoint in
the behavior of interpolating sequences for $\hf$.

We shall now turn to the case of big Hardy-Orlicz spaces.
%We can now affirm that this behavior in fact changes between exponential
%and polynomial growth. 
Let $\vp$ be a strongly convex function
with associated Hardy-Orlicz space $\hf$. 
Our central assumption on
$\vp$ is the quasi-triangular inequality
\bea\label{deltaineq}
      \vp(a+b)\le c(\vp(a)+\vp(b))
\eea
for some fixed constant $c\ge 1$ and for all $a,b\ge t_0$.
This condition obviously implies that $\hf$ is stable with respect
to multiplication so that under this condition $\hf$ is an algebra.
Observe that \eqref{deltaineq} is an equivalent formulation 
of the usual $\Delta_2$-condition \eqref{Delta} (now for 
$\vp$ instead of $\Phi$), 
and we will henceforth denote the condition
%and this will be our notation for
%to which we will refer for 
% will be the notation for 
\eqref{deltaineq}
%and we shall thus use the notation
%$\Delta_2$ %refer to $\Delta_2$-condition
by $\Delta_2$.
%throughout the rest of the paper.
Note that the $\Delta_2$-condition obviously implies $\tilde{\Delta}_2$.
%for $\vp$.

%The standard example in this setting is $\vp_p(t)=t^p$ for $p>1$
%and e.g.\ $t\ge 1$
%(the equation \eqref{deltaineq} holds also for $p\in (0,1]$ but
%this range is excluded by the $\Delta_2$-condition).
%The condition \eqref{deltaineq} entails that $\vp$ grows at most polynomially.

%For the condition that gives the description of the existence
%of a harmonic majorant in $L^{\vp}$
%w
%We well need a
Another condition on $\vp$ will be used. We say that $\vp$ satisfies
the $\nabla_2$-condition (see e.g.\ \cite{les}) if there exist
$d>1$ and $t_0>0$ such that
\beqa
 2\vp(t)\le \frac{1}{d}\vp(dt),\quad t\ge t_0,
\eeqa
(see Section \ref{HO} for more on this).

The main result of this paper then reads as follows.

\begin{theorem}\label{CNS}
Let $\vp:\R\lra [0,\infty)$ be a strongly convex function 
satisfying the $\Delta_2$-condition. The following assertions are
equivalent.
\begin{itemize}
\item[(a)] $\Lambda$ is a free interpolating sequence
for $\hf$: $\Lambda\in \Int \hf$
\item[(b)] There exists a positive measurable function
%weight $w\in L^1(\T)$ 
%such that
%
$w\in L^{\vp} (\T)$ such that $\varphi_\Lambda\le P[w]$.
\item[(c)] The trace space is given by
\beqa
 \hf|\Lambda=l_{\Phi}:=\{(a_{\lambda}):\exists 0\le w\in L^{\vp}(\T)
 \text{ %measurable
   with } %\vp\circ w \in L^{1},
 \log^+|a_{\lambda}|\le P[w](\lambda)\}.
\eeqa
\end{itemize}

If moreover $\vp$ satisfies the $\nabla_2$-condition then the above
three conditions are equivalent to the following.

\begin{itemize}
\item[(d)] There exists a constant $C>0$ such that for any sequence of
non-negative numbers $(c_{\lambda})$,
\beqa
 \sum_{\lambda\in\Lambda}c_{\lambda}\log\frac{1}{|B_{\lambda}(\lambda)|}
 \le C
 \|\sum_{\lambda\in\Lambda}c_{\lambda}P_{\lambda}(\zeta)\|_{(L^{\vp})^*}
\eeqa
\end{itemize}
\end{theorem}

The standard examples of functions satisfying $\Delta_2$ and $\nabla_2$ are
$\vp(t)=\vp_p(t):=t^p$ for fixed $p>1$ and $t\ge t_0$,
or $\vp=\psi_{\eps}(t):=t\log^{\eps}t$ 
for some fixed $\eps>0$ and $t\ge t_0$. Note that $\vp_p$ satisfies
both conditions also for $p\in (0,1]$ but this range is excluded by
the definition of strongly convex functions. 

The space $L^{\vp}$ appearing in the theorem is the standard Orlicz space
of measurable functions $u$ such that $\vp\circ |u|\in L^1(\T)$.
As a consequence of the $\Delta_2$-condition \eqref{deltaineq}
it turns out that its dual space
$(L^{\vp})^*$ 
% is its dual,
%which, as it turns out,  % that by the $\delta_2$ condition, the space
%$(L^{\vp})^*$ 
is in fact also an Orlicz space (associated
with the complementary function of $\vp$, see Section \ref{HO}
for more comments and details).

It is interesting to note the analogy between condition (d) of
this theorem and \eqref{charNevInt} which characterizes the
interpolating sequences for the Nevanlinna class. Recall that 
\eqref{charNevInt}
was not sufficient for the existence of a quasi-bounded
harmonic majorant. In our situation however, 
%as soon as
%the defining function has 
any growth strictly faster than in the
$L^1$-situation %, the condition (d)
suffices to eliminate the singular part of the measure defining the
harmonic majorant (this is maybe not so surprising, there is a kind
of ``de la Vall\'ee Poussin effect'', see \cite[Theorem 4.14]{RR}).
%to give a quasi-bounded majorant (and 
%(which is in the right Orlicz space).

As in \cite{HMNT} we will investigate the problem of existence
of harmonic majorants in the general setting. More precisely we
are interested in the question
when a Borel function defined on $\DD$ admits a
harmonic majorant $P[w]$ with $w\in L^{\vp}$. The answer to this
problem is given by the following result which involves
the so-called Poisson balayage. Recall that the Poisson balayage of a finite
positive measure $\mu$ in the closed unit disk is defined as
\beqa
 B(\mu)(\zeta)=\int_{\DD}P_z(\zeta)\,d\mu(z),\quad \zeta\in\T.
\eeqa

\begin{theorem}\label{charharmmaj}
Let $\vp$ be a strongly convex function 
%such that \eqref{deltaineq} holds and 
that satisfies the $\Delta_2$-condition and the $\nabla_2$-condition.
If $u$ is a non-negative Borel function on the unit disk then
the following two assertions are equivalent.

\begin{itemize}
\item[(a)] There exists a function $w\in L^{\vp}$ such that
$u(z)\le P[w](z)$ for all $z\in \DD$.
\item[(b)] There exists a constant $C\ge 0$ such that
\beqa
 \sup_{\mu\in {\mathcal B}_{\vp^*}} \int u(z)\,d\mu(z)\le C,
\eeqa
where ${\mathcal B}_{\vp^*}=\{\mu $ positive measure on $\DD$: 
$\|B\mu\|_{(L^{\vp})^*}\le 1\}$.
\end{itemize}
\end{theorem}

It is again interesting to point out the analogy between this result
and that given in \cite[Theorem 1.4]{HMNT}. %However
Note also that the corresponding condition (b) in that theorem
does not give a quasi-bounded majorant but only a
harmonic majorant, and a more subtle 
condition is needed to handle the case of quasi-bounded 
majorants (see \cite[Theorem 1.6]{HMNT}).

The paper is organized as follows.
In the next section we shall add some more comments on 
Orlicz and Hardy-Orlicz spaces.
The sufficiency part of our main theorem has been given in
\cite[Theorem 9.1]{HMNT} and we refer the reader to that
paper for a proof. 
The general structure of the proof of the necessity
goes along the lines of the necessary part for the Smirnov class.
However, the key result \cite[Proposition 4.2]{HMNT} does not work
any longer in our context. Note that it is precisely that proposition 
which shows that separated Blaschke
sequences are interpolating for the Nevanlinna and
Smirnov classes. We will actually discuss in some details
in Section \ref{examples} an example showing that 
such a result cannot be expected in Hardy-Orlicz classes. 
The example shows that there are big Hardy-Orlicz 
spaces which are close to each other
in a sense and for which there exist even separated Blaschke
sequences that are interpolating for one space but not for the other one.
For this reason we need a new idea %at this point 
which will be discussed in Section \ref{necessity}.
The key is to factorize the Blaschke product $B_{\Lambda}$ into two factors 
that behave essentially in the same way as $B_{\Lambda}$
(in a sense to be made precise). This will be
achieved through a theorem by Hoffman.
The trace space characterization is quite immediate and will be
discussed in Section \ref{trace}.
Concerning harmonic majorants, we will discuss this problem
in Section \ref{harmmaj}. We have already insisted in the analogy
between our Theorem \ref{charharmmaj} and Theorems 1.4 and 1.6 in
\cite{HMNT}. The techniques that apply in our situation are
more classical than those used in \cite{HMNT}: we will use some
duality arguments and a theorem by Mazur-Orlicz on positive
linear functionals (which is essentially the Hahn-Banach theorem).
The equivalence of (b) and (d) of Theorem \ref{CNS} 
then follows from Theorem \ref{charharmmaj} (it suffices to consider
positive measures $\mu$ supported on $\Lambda$).

\begin{questions*}
With our big Hardy-Orlicz spaces we narrow down the gap in the
description of interpolating sequences from
above: coming from the Smirnov class where the harmonic majorant 
must have an absolutely continuous measure, and so with an integrable
weight $w\in L^1$, we obtain now that the weight has to be in
$L^{\vp}$ (and we can get in a way arbitrarily close to $L^1$, see
also Section \ref{examples}). Now, we have already indicated that the
$\Delta_2$-condition \eqref{deltaineq} implies a polynomial growth on $\vp$.
Observe that for classical Hardy spaces $H^p$ --- where the Carleson
condition characterizes the interpolating sequences --- the
defining functions are given by $t\lmto e^{pt}$. So one could ask
what happens for defining functions with subexponential growth
like e.g.\ $\vp(t)=e^{\sqrt{t}}$ and whether there are Hardy-Orlicz
spaces beyond $H^p$, $p>0$, for which the Carleson condition still
characterizes the interpolating sequences.
\end{questions*}

\section{(Hardy-) Orlicz spaces}\label{HO}

For the background on the notions related to Hardy-Orlicz classes
used in this section
we refer to \cite{RR}, \cite{les} and \cite{KR}.

%The 
Let $\vp$ be a strongly convex function.
One way of introducing the Hardy-Orlicz class $\hf$ is to
take all the functions $f\in N^+$ such that the subharmonic function
$\vp(\log^+|f|)$ admits a harmonic majorant on $\DD$
(see \cite[Definition 3.15]{RR}). 
We have already introduced the notation
\beqa
 \Phi=\vp\circ\log
\eeqa
keeping in mind that the function $\Phi$ is chosen in a way
guaranteeing that $\Phi(|f|)$ is subharmonic.

Since $f$ is in the Smirnov class,
the fact that $\Phi(|f|)$ has a harmonic majorant is equivalent to
(see \cite[Theorem 4.18]{RR})
\beqa
 J_{\Phi}(f)=\int_{\T}\Phi(|f(e^{it})|)\;dt
 =\int_{\T}\vp(\log|f(e^{it})|)\;dt<\infty,
\eeqa
so that $\hf$ can be defined as
\beqa
 \hf=\{f\in N^+:J_{\Phi}(f)=\int_{\T}\vp(\log|f(e^{it})|)\;dt<\infty\}.
\eeqa
Observe that $\hf$ does not depend on the behaviour of $\vp$ 
for $t\le t_0$ whenever a $t_0\in\R$ is fixed.

The integral expression $J_{\Phi}$ is called a modular, and it does not 
define a metric on $\hf$. A metric can be defined by $d(f,g)=
\|f-g\|_{\Phi}$ where
\bea\label{metric}
 \|f\|_{\Phi}=\inf\{t>0:J_{\Phi}(f/t)\le t\},
\eea
and $\hf$ equipped with this metric is a complete space.

In our situation, thanks to the $\Delta_2$-condition, 
$J_{\Phi}(af)<\infty$ for any $a>0$ when
$J_{\Phi}(f)<\infty$ so that we do not need to distinguish
between the Orlicz class, the Orlicz space and the space of 
finite elements (in the terminologie of \cite{les}).

It is clear that if $f\in \hf$ then by the Riesz-Smirnov
factorization $f=Ih$ where $I$ is an inner function and
$h$ is outer in $N^+$. Clearly $J_{\Phi}(f)=J_{\Phi}(h)$.
Since $h$ is outer in $N^+$ we have
$h(z)=\exp(\int (\zeta+z)/(\zeta-z) w(\zeta)\,dm(\zeta))$ for some
real function $w\in L^1(\T)$, so that $|h|=\exp(P[w])$ in $\DD$ which
has boundary values $\exp(w)$ $m$-almost everywhere on $\T$. Hence
\beqa
 J_{\Phi}(h)=\int_{\T}\vp(\log(\exp w))\;dm=\int_{\T} \vp(w)\;dm
 \ge \int_{\T} \vp(w_+)\;dm=J_{\vp}(w_+),
\eeqa
where $w_+=\max(0,w)$.
In other words $w_+$ is in the Orlicz class $L^{\vp}$ of measurable
functions $u$ such that $J_{\vp}(u)<\infty$. 
%Note that the function
%$\vp$ is convex increasing so that $L^{\vp}$ is a Banach space.

Moreover $\lim_{t\to\infty}\vp(t)/t=+\infty$, and $L^1=L^{\vp_1}$
(recall that $\vp_1(t)=t$, $t\ge 1$),
so that by standard results on Orlicz spaces
we get $L^{\vp}\subset L^1$. In particular,
if we now take any $w\in L^{\vp}$, then $w\in L^1$, and we can define an
outer function in the Smirnov class by 
\beqa
 f_w(z)=\exp\left(\int
   \frac{\zeta+z}{\zeta-z}w(\zeta)\;dm(\zeta)\right).
\eeqa
%which is of course in the Smirnov class.
For the same reasons as above this function has boundary limits 
$\exp(w)$ a.e.\ on $\T$ and 
\beqa
 J_{\Phi}(f_w)&=&\int_{\T}\vp(w)\;dm=\int_{w\ge 0} \vp(|w|)\;dm
 +\int_{w<0} \vp(w)\;dm\\
 &=&\int \vp(|w|)dm+\int_{w<0}\vp(w)-\vp(|w|)\;dm
 \le \vp(0)+ J_{\vp}(w)<\infty,
\eeqa 
so that
$w$ gives rise to an outer function in $\hf$.

Another important fact on Hardy-Orlicz classes that will be useful
for us later is an estimate on point evaluations.

Indeed, the maximal radial growth that we can attain for $f\in\hf$ is
\beqa
 |f(z)|\le \Phi^{(-1)}\left(\frac{J_{\Phi}(f)}{1-|z|}\right).
\eeqa
This can be deduced from the subharmonicity of $\Phi(|f|)$ and
the change of variable $u\lmto b_z(u)$ (see also \cite[II.1.2]{les}).
Hence %if $f\in \Ball_{\hf}$ then
\bea\label{pointeval}
 \log |f(z)|\le \vp^{-1}\left(\frac{c_f}{1-|z|}\right).
\eea
For classical Hardy spaces $H^p$ one recovers from this the
usual estimate on the point evaluation. 

We will use the above estimate to show that certain separate
sequences are not interpolating for ``big'' Hardy-Orlicz classes (see Section
\ref{examples}).

Some more tools need to be introduced in the context of Orlicz spaces.
In the above arguments, we did not need to appeal to topological
considerations in $L^{\vp}$. Thus the definition of $L^{\vp}$ 
as the set of measurable functions $w$ for which $\vp\circ |w|\in L^1$
was sufficient. Later on however
we will have to consider in particular the dual and
the bidual of $L^{\vp}$, and hence we need a proper norm in $L^{\vp}$. 
The alerted reader might have observed that $\vp$ is not a
so-called $N$-function since for instance it does not vanish at $0$ and hence
it is not appropriated to define a topology on $L^{\vp}$. 
In the sequel, when considering $L^{\vp}$ we will
think of $\vp$ as being suitably replaced on $[0, t_0]$ in order that
the resulting function regarded as a function on $[0,+\infty)$
is convex and vanishing conveniently in $0$. Then, in view of the convexity
of $\vp$ the Orlicz space $L^{\vp}$ equipped with the metric
\eqref{metric} ($\Phi$ replaced by $\vp$) is a Banach space.

%Recall that our function $\vp$ is supposed to be convex. Hence it
%admits (e.g.) a right derivative $p$ in every point which is right
%continuous. Setting $q(s)=\sup_{p(t)\le s}t$, and setting

With the adjusted function $\vp$ we can define the so-called
complementary function. It is defined by
$\vp^*(s)=\max_{t\ge 0}\{st-\vp(t)\}$ (it is also possible to
define $\vp^*$ using the ``inverse'' of the right derivative
of $\vp$).
%$\vp^*(t)=\int_0^t q(s)\,ds$ one gets the so-called complementary
%function of $\vp$. 
Since $\vp$ satisfies the $\Delta_2$-condition, 
we get $(L^{\vp})^*=L^{\vp^*}$ (see for instance \cite{KR} for this).
Note also that $(\vp^*)^*=\vp$.

%The function $\vp$ is said to satisfy
We have also mentioned the $\nabla_2$-condition. Recall that
$\vp$ satisfies the $\nabla_2$-condition
if there are constants $d>1$ and $t_0\ge 0$ such that
for all $t\ge t_0$ we have $2\vp(t)\le \vp(dt)/d$. This condition
is actually equivalent to the fact that $\vp^*$ satisfies the
$\Delta_2$-condition (\cite[Chapter 1, Theorem 4.2]{KR}) so that
$(L^{\vp})^{**}=(L^{\vp^*})^*=L^{\vp^{**}}=L^{\vp}$. 
%since moreover $\vp^{**}=\vp$. 
In particular, if the strongly 
convex function $\vp$ satisfies both the $\Delta_2$-condition
and the $\nabla_2$-condition, then
$L^{\vp}$ is a reflexive space (and the converse is also true,
see \cite[p.56]{les}).

%The alerted reader has probably observed that $\vp$ is not a
%so-called $N$-function since it does not vanish in $0$ and hence
%it is not appropriated to define a topology on $L^{\vp}$. 
%However, it is possible to replace $\vp$ on $[0, t_0]$ without
%changing the set $L^{\vp}$ and such that the corresponding function is an
%$N$-function. The above reasoning applies then to the modified
%function $\vp$. It should also be noted that the complementary
%function can be given by $\vp^*(s)=\max_{t\ge 0}\{st-\vp(t)\}$, $s\ge 0$.

%{\sc add something on dual Orlicz spaces, complementary functions}

\section{Necessary condition}\label{necessity}

The central Proposition 4.1 of \cite{HMNT} which claims that
$\log (1/|B(z)|)$ is controlled by a (quasi-bounded) positive
harmonic function for those $z$ which are uniformly bounded away in 
the pseudohyperbolic metric from the zeros of $B$ cannot be applied
to our situation since the weight of the quasi-bounded majorant
need not be in $L^{\vp}$. In 
Section \ref{examples} we will give examples
of separated sequences for which
the majorants are in no $L^{\vp}$ whenever $\vp(t)\ge t\log^{\eps} t$ for
$t\ge t_0$ and %some positive constant 
$\eps>0$.
%We thus need a more delicate
%tool to come over this difficulty. 
To overcome this difficulty we need a more precise results which
is based on a theorem by Hoffman
%
%In order to prove the necessary condition we will check
%that apart from Proposition 4.1 the techniques used in \cite{HMNT} apply
%to our situation. It is essentially Hoffman's theorem 
%that we shall state here 
(see \cite[p. 411]{Gar}):

\begin{theorem}[Hoffman's theorem]\label{hoffthm}
For $0<\delta<1$ there are constants $a=a(\delta)$ and $b=b(\delta)$
such that the Blaschke product $B(z)$ with zero set
$\Lambda$ has a nontrivial factorization
$B=B_1B_2$ such that
\beqa
 a|B_1(z)|^{1/b}\le |B_2(z)|\le \frac{1}{a}|B_1(z)|^b
\eeqa
for every $z\in \DD\setminus \bigcup_{\lambda\in\Lambda}D(\lambda,\delta)$
where $D(\lambda,\delta)=\{z\in\DD:|b_{\lambda}(z)|<\delta\}$ 
is the pseudohyperbolic disk centered at $\lambda$ with radius $\delta$.
\end{theorem}

\begin{corollary}\label{splitting}
Let $\Lambda=\{\lambda_n\}_n\subset\DD$ be a separated Blaschke
sequence. Then there exists a partition 
\beqa
 \Lambda=\Lambda_1\stackrel{\cdot}{\cup}\Lambda_2
\eeqa
and constants $c,\eta>0$ such that
\bea\label{estim}
 \log\frac{1}{|(B_k)_{\lambda}(\lambda)|}
 \ge c \log\frac{1}{|B_{\lambda}(\lambda)|}-\eta,
\eea
where $B_k=B_{\Lambda_k}=\prod_{\mu\in \Lambda_k}b_{\mu}$ and
$(B_k)_{\lambda}=B_{\Lambda_k\setminus \{\lambda\}}$ if $\lambda\in
\Lambda_k$, $(B_k)_{\lambda}=B_k$ otherwise.
\end{corollary}

\begin{proof}
One should first note that the constants in Hoffman's theorem
depend only on $\delta$. 

Let $\delta$ be the separation constant of the sequence $\Lambda$ and
fix $\Lambda=\Lambda_1\stackrel{\cdot}{\cup}\Lambda_2$ a partition of 
$\Lambda$ obtained from Hoffman's theorem.

Pick $\lambda\in \Lambda$, say $\lambda\in\Lambda_1$ (so that in the
following considerations we will assume $k=1$).
%According to Hoffman's theorem the sequence
Then
$\Lambda\setminus\{\lambda\}=\Lambda_1\setminus\{\lambda\}\cup\Lambda_2$.
A careful inspection of the proof of Hoffman's theorem (see e.g.\ the
indicated reference) shows that we have
\bea\label{hoff}
 a|B_{\Lambda_1\setminus\{\lambda\}}(z)|^{1/b}\le 
 |B_{\Lambda_2}(z)|\le \frac{1}{a}|B_{\Lambda_1\setminus\{\lambda\}}(z)|^b
\eea
%with the same constants as for the initial splitting $B=B_1B_2$, but
%now the estimates are valid 
for $z\in \DD\setminus\bigcup_{\mu\in
\Lambda\setminus\{\lambda\}}D(\mu,\delta)$, and so in particular 
for $\lambda$. 

We should pause here to add some comments on the proof of Hoffman's
theorem given in \cite{Gar}. % for our situation. 
The fact that we take away
$\lambda$ implies a possible shift between the odd and the even
indexed points in the strip $T_k$ (according to the terminology
in \cite{Gar}) containing $\lambda$, and so 
the choice of $\Lambda_1$ and $\Lambda_2$ depends on $\lambda$.
However, this is of no harm since these shifts mean in fact that we
just add or take away at most one term in each layer $T_k$. 
Since the layers $T_k$ are of pseudohyperbolic
constant thickness, the terms added or subtracted correspond to
an $\Hi$-interpolating sequence the constant
(in the sense of \eqref{Intcst}) of which is bounded
by that of an $\Hi$-interpolating sequence in a radius with given
separation constant. So in the estimates \eqref{hoff}, the constant
$b$ is the same as in the splitting of the original sequence $\Lambda=
\Lambda_1\stackrel{\cdot}{\cup}\Lambda_2$ 
whereas the constant $a$ should be replaced
by a different one, but independent on $\lambda$.

So
\beqa
 \log\frac{1}{|B_{\Lambda\setminus\{\lambda\}}(\lambda)|}
 &=&\log\frac{1}{|B_{\Lambda_1\setminus\{\lambda\}}(\lambda)|}+
  \log\frac{1}{|B_{\Lambda_2}(\lambda)|}\\
 &\le& \log\frac{1}{|B_{\Lambda_1\setminus\{\lambda\}}(\lambda)|}
 +\log\frac{1}{a|B_{\Lambda_1\setminus\{\lambda\}}(\lambda)|^{1/b}}\\
%
%\log\frac{1}{a}+\frac{1}{|B_{\Lambda_2}(\lambda)|}+
% \log\frac{1}{|B_{\Lambda_1\setminus\{\lambda\}}(\lambda)|}\\
 &=&\frac{b+1}{b}
 \log\frac{1}{|B_{\Lambda_1\setminus\{\lambda\}}(\lambda)|}-\log{a},
\eeqa
% \log\frac{1}{|B_{\Lambda_1\setminus\{\lambda\}}(\lambda)|^{b+1}}
and we can set $c=(b+1)/b$ and $\eta=\log a$.

The cases $\lambda\in\Lambda_2$, $k=2$ are treated in a similar way.
\end{proof}

We are now in a position to state the desired result for separated sequences.

\begin{corollary}\label{intsep}
If $\Lambda$ is a separated sequence that is interpolating for
$\hf$ then 
%there is a quasibounded harmonic majorant for
%$\vp_{\Lambda}$ in $L^{\vp}$, i.e.\ 
there is a positive measurable
function $w$ with $\vp\circ w\in L^{1}$ such
that
\beqa
 \log\frac{1}{|B_{\lambda}(\lambda)|}\le P[w](\lambda),\quad \lambda\in
 \Lambda.
\eeqa
\end{corollary}

%Using the fact that $h:=\exp(\frac{1}{2\pi}\int\frac{\zeta+z}{\zeta-z}
%(-w(\zeta))dm(\zeta))\in \hf$, we can give
%another formulation of the result of the Theorem, which
%will be useful later on.
%
%{\it If $\Lambda$ is a separated interpolating sequence for
%$\hf$ then there existe an outer function $h\in \hf$ such that
%$|B_{\lambda}(\lambda)|\ge |h(\lambda)|$}.

\begin{proof}
Suppose that $\Lambda$ is a separated sequence.
% interpolating sequence for $\hf$. 
By Corollary \ref{splitting}, there exists a partition
\beqa
 \Lambda=\Lambda_1\stackrel{\cdot}{\cup}\Lambda_2,\quad \lambda\in\Lambda,
\eeqa
and constants $c,\eta>0$
such that for all $\lambda\in \Lambda$
%\beqa
\bea\label{estim1}
 \log\frac{1}{|(B_k)_{\lambda}(\lambda)|}
 \ge c \log\frac{1}{|B_{\lambda}(\lambda)|}-\eta.
\eea
%where $B_k=B_{\Lambda_k}=\prod_{\mu\in \Lambda_k}b_{\mu}$, $k=1,2$.
Since $\Lambda$ is moreover interpolating for $\hf$ there 
exist two functions $f_i\in \hf$, $i=1,2$, such that
\beqa
 f_i|\Lambda_i=1,\quad f_i|(\Lambda\setminus\Lambda_i)=0.
\eeqa
Now $\hf\subset N^+$, and we can factorize in the following way
\beqa
 f_i=B_{\Lambda\setminus\Lambda_i}I_ih_i, \quad i=1,2,
\eeqa
where $I_i$ is an inner function, $h_i$ is outer in $\hf$:
\beqa
 h_i(z)=\exp\left(\frac{1}{2\pi}\int\frac{\zeta+z}{\zeta-z}
 w_i(\zeta)dm(\zeta)\right),
\eeqa
and $(w_i)_+\in L^{\vp}$. %Clearly $w_i$ can be supposed nonnegative.
Then
for every $\lambda\in \Lambda_i$
\beqa
 1=f_i(\lambda)=|f_i(\lambda)|\le |B_{\Lambda\setminus\Lambda_i}(\lambda)|
 \cdot |h_i(\lambda)|,
\eeqa
so that
\beqa
 \log\frac{1}{|B_{\Lambda\setminus\Lambda_i}(\lambda)|}
 \le P[w_i](\lambda)\le P[(w_i)_+](\lambda).
\eeqa
Using \eqref{estim1}, we get
\beqa
  \log\frac{1}{|B_{\lambda}(\lambda)|}
 \le \frac{1}{c}\left(\eta+
 \log\frac{1}{|B_{\Lambda\setminus\Lambda_i}(\lambda)|}
 \right)
 \le P[\frac{1}{c}((w_i)_++\eta)](\lambda).
\eeqa
The corollary then follows by setting $w=\frac{1}{c}((w_1)_++
(w_2)_++\eta)$ which is still in $L^{\vp}$.
\end{proof}

Let us now switch to the necessary condition in the general
situation.
The central trick in \cite{HMNT} is to decompose the disk $\DD$ into 
Whitney ``cubes'' that split
the sequence into four pieces that are uniformly separated from
each other in the pseudohyperbolic metric. This allows one to
reduce the situation to the separated one. 

We will repeat here the proof given in \cite{HMNT} for completeness
adding the necessary changes for the situation of big Hardy-Orlicz spaces. 
\\
\\
%\begin{center}
%\input{dyadicpartition.latex}\\
%Figure: dyadic partition
%\end{center}
\begin{center}
\begin{picture}(0,0)%
\includegraphics{dyadicpart.pstex}%
\end{picture}%
\setlength{\unitlength}{1657sp}%
\begingroup\makeatletter\ifx\SetFigFont\undefined%
\gdef\SetFigFont#1#2#3#4#5{%
  \reset@font\fontsize{#1}{#2pt}%
  \fontfamily{#3}\fontseries{#4}\fontshape{#5}%
  \selectfont}%
\fi\endgroup%
\begin{picture}(7764,7224)(2644,-9298)
\put(4006,-3796){\makebox(0,0)[lb]{\smash{{\SetFigFont{5}{6.0}{\familydefault}{\mddefault}{\updefault}{\color[rgb]{0,0,0}1}%
}}}}
\put(4006,-6496){\makebox(0,0)[lb]{\smash{{\SetFigFont{5}{6.0}{\familydefault}{\mddefault}{\updefault}{\color[rgb]{0,0,0}3}%
}}}}
\put(7471,-6451){\makebox(0,0)[lb]{\smash{{\SetFigFont{5}{6.0}{\familydefault}{\mddefault}{\updefault}{\color[rgb]{0,0,0}4}%
}}}}
\put(3376,-7981){\makebox(0,0)[lb]{\smash{{\SetFigFont{5}{6.0}{\familydefault}{\mddefault}{\updefault}{\color[rgb]{0,0,0}1}%
}}}}
\put(5266,-7981){\makebox(0,0)[lb]{\smash{{\SetFigFont{5}{6.0}{\familydefault}{\mddefault}{\updefault}{\color[rgb]{0,0,0}2}%
}}}}
\put(7066,-7981){\makebox(0,0)[lb]{\smash{{\SetFigFont{5}{6.0}{\familydefault}{\mddefault}{\updefault}{\color[rgb]{0,0,0}1}%
}}}}
\put(8641,-7936){\makebox(0,0)[lb]{\smash{{\SetFigFont{5}{6.0}{\familydefault}{\mddefault}{\updefault}{\color[rgb]{0,0,0}2}%
}}}}
\put(3241,-8656){\makebox(0,0)[lb]{\smash{{\SetFigFont{5}{6.0}{\familydefault}{\mddefault}{\updefault}{\color[rgb]{0,0,0}3}%
}}}}
\put(4141,-8656){\makebox(0,0)[lb]{\smash{{\SetFigFont{5}{6.0}{\familydefault}{\mddefault}{\updefault}{\color[rgb]{0,0,0}4}%
}}}}
\put(5041,-8656){\makebox(0,0)[lb]{\smash{{\SetFigFont{5}{6.0}{\familydefault}{\mddefault}{\updefault}{\color[rgb]{0,0,0}3}%
}}}}
\put(5941,-8656){\makebox(0,0)[lb]{\smash{{\SetFigFont{5}{6.0}{\familydefault}{\mddefault}{\updefault}{\color[rgb]{0,0,0}4}%
}}}}
\put(6841,-8656){\makebox(0,0)[lb]{\smash{{\SetFigFont{5}{6.0}{\familydefault}{\mddefault}{\updefault}{\color[rgb]{0,0,0}3}%
}}}}
\put(7696,-8656){\makebox(0,0)[lb]{\smash{{\SetFigFont{5}{6.0}{\familydefault}{\mddefault}{\updefault}{\color[rgb]{0,0,0}4}%
}}}}
\put(8641,-8656){\makebox(0,0)[lb]{\smash{{\SetFigFont{5}{6.0}{\familydefault}{\mddefault}{\updefault}{\color[rgb]{0,0,0}3}%
}}}}
\put(9541,-8656){\makebox(0,0)[lb]{\smash{{\SetFigFont{5}{6.0}{\familydefault}{\mddefault}{\updefault}{\color[rgb]{0,0,0}4}%
}}}}
\put(7201,-9016){\makebox(0,0)[lb]{\smash{{\SetFigFont{5}{6.0}{\familydefault}{\mddefault}{\updefault}{\color[rgb]{0,0,0}2}%
}}}}
\put(7651,-9016){\makebox(0,0)[lb]{\smash{{\SetFigFont{5}{6.0}{\familydefault}{\mddefault}{\updefault}{\color[rgb]{0,0,0}1}%
}}}}
\put(8101,-9016){\makebox(0,0)[lb]{\smash{{\SetFigFont{5}{6.0}{\familydefault}{\mddefault}{\updefault}{\color[rgb]{0,0,0}2}%
}}}}
\put(8551,-9016){\makebox(0,0)[lb]{\smash{{\SetFigFont{5}{6.0}{\familydefault}{\mddefault}{\updefault}{\color[rgb]{0,0,0}1}%
}}}}
\put(9001,-9016){\makebox(0,0)[lb]{\smash{{\SetFigFont{5}{6.0}{\familydefault}{\mddefault}{\updefault}{\color[rgb]{0,0,0}2}%
}}}}
\put(9451,-9016){\makebox(0,0)[lb]{\smash{{\SetFigFont{5}{6.0}{\familydefault}{\mddefault}{\updefault}{\color[rgb]{0,0,0}1}%
}}}}
\put(9901,-9016){\makebox(0,0)[lb]{\smash{{\SetFigFont{5}{6.0}{\familydefault}{\mddefault}{\updefault}{\color[rgb]{0,0,0}2}%
}}}}
\put(6751,-9016){\makebox(0,0)[lb]{\smash{{\SetFigFont{5}{6.0}{\familydefault}{\mddefault}{\updefault}{\color[rgb]{0,0,0}1}%
}}}}
\put(6301,-9016){\makebox(0,0)[lb]{\smash{{\SetFigFont{5}{6.0}{\familydefault}{\mddefault}{\updefault}{\color[rgb]{0,0,0}2}%
}}}}
\put(5851,-9016){\makebox(0,0)[lb]{\smash{{\SetFigFont{5}{6.0}{\familydefault}{\mddefault}{\updefault}{\color[rgb]{0,0,0}1}%
}}}}
\put(5356,-9016){\makebox(0,0)[lb]{\smash{{\SetFigFont{5}{6.0}{\familydefault}{\mddefault}{\updefault}{\color[rgb]{0,0,0}2}%
}}}}
\put(4906,-9016){\makebox(0,0)[lb]{\smash{{\SetFigFont{5}{6.0}{\familydefault}{\mddefault}{\updefault}{\color[rgb]{0,0,0}1}%
}}}}
\put(4456,-9016){\makebox(0,0)[lb]{\smash{{\SetFigFont{5}{6.0}{\familydefault}{\mddefault}{\updefault}{\color[rgb]{0,0,0}2}%
}}}}
\put(4006,-9016){\makebox(0,0)[lb]{\smash{{\SetFigFont{5}{6.0}{\familydefault}{\mddefault}{\updefault}{\color[rgb]{0,0,0}1}%
}}}}
\put(3556,-9016){\makebox(0,0)[lb]{\smash{{\SetFigFont{5}{6.0}{\familydefault}{\mddefault}{\updefault}{\color[rgb]{0,0,0}2}%
}}}}
\put(3106,-9016){\makebox(0,0)[lb]{\smash{{\SetFigFont{5}{6.0}{\familydefault}{\mddefault}{\updefault}{\color[rgb]{0,0,0}1}%
}}}}
\end{picture}%
\\
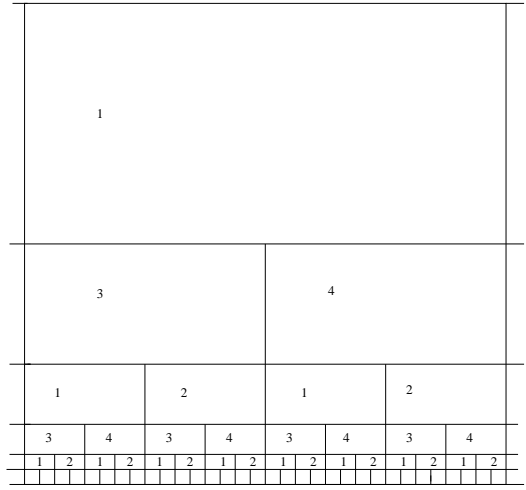
Figure 1: dyadic partition
\end{center}

\vspace{0.5cm}

Let $I_{n,k}:=\{ \eit : \theta \in [{2\pi}k 2^{-n},{2\pi}(k+1)
 2^{-n}) \}, \ 0 \le k < 2^n$, be the dyadic arcs
and $Q_{n,k}:= \{ r \eit : \eit \in I_{n,k} , 1- 2^{-n} \le r < 1-
 2^{-n-1} \}$ the associated  ``dyadic squares".

%Observe that the hyperbolic diameter of each Whitney square
%$Q_{n,k}$ is bounded between two absolute constants.

The splitting of the sequence into four pieces
will be done in the following way: $\Lambda=
\bigcup_{i=1}^4 \Lambda_i$ such that each piece $\Lambda_i$ lies in a union
of dyadic squares
% as in \eqref{dyadsquares}
that are uniformly separated from each other
(see Figure 1).
More precisely, set
\beqa
       \Lambda_1=\Lambda\cap Q^{(1)},
\eeqa
where the family $Q^{(1)}$ is given by $\{Q_{2n,2k}\}_{n,k}$
(for the remaining three sequences we respectively choose
$\{Q_{2n ,2k+1}\}_{n,k}$,
$\{Q_{2n+1,2k}\}_{n,k}$ and
$\{Q_{2n+1,2k+1}\}_{n,k}$).
In order to avoid technical difficulties we count only those
$Q$ containing points of $\Lambda$.
% (in case $\Lambda \cap Q$
%is empty there is nothing to prove).
In what follows we will argue on one sequence, say $\Lambda_1$.
The arguments are the same for the other sequences.

%Our first observation is that,
By construction,
for $Q, L \in Q^{(1)}$, $Q \neq L$,
\beqa
       \rho(Q,L):=\inf_{z\in  Q, w\in L}\rho(z,w)
       \ge\delta>0,
\eeqa
for some fixed $\delta$. In what follows, the letters $j$, $k$... will
stand for indices in $\mathbb N^2$ of the form $(n,l), 0 \le l <2^n$.
The closed rectangles
$\overline{Q_j}$ are compact in
$\DD$ so that $\Lambda_1\cap Q_j$ can only
contain a finite number of points (they contain at least one point,
by assumption).
Therefore
\beqa
       0<m_j :=\min_{\lambda\in \Lambda_1\cap Q_j} |B_{\lambda}(\lambda)|
\eeqa
(note that we consider the entire Blaschke product $B_\lambda$
associated with $\Lambda\setminus\{\lambda\}$). Take
$\lambda_j^{1}\in Q_j$ such that
$m_j=|B_{\lambda_j^{1}}(\lambda_j^{1})|$.

Assume now that $\Lambda\in \Int \hf$.
So, since $\loi\subset \hf\vert\Lambda$, there exists a function $f_1\in \hf$
such that
\beqa
       f_1(\lambda)=
       \begin{cases}
1\quad &\text{if}\quad \lambda\in \{\lambda_j^{1}\}_j\\
0\quad &\text{if}\quad \lambda\in
\Lambda\setminus\{\lambda_j^{1}\}_j.
\end{cases}
\eeqa
By the Riesz-Smirnov factorization we have
\begin{equation}
\label{factorizn}
       f_1=B_{\Lambda\setminus \{\lambda_j^{1}\}_j}I_1g_1,
\end{equation}
where $I_1$ is some inner function and $g_1$ is outer in $\hf$.
Hence there exists a weight $w_1\in L^1$ such that $g_1=f_{w_1}$
and $(w_1)_+\in L^{\vp}$.
As in the proof of Corollary \ref{intsep} we get
\beqa
       1 = |f_1(\lambda_k^{1})|
        \le |B_{\Lambda\setminus \{\lambda_j^{1}\}_j}(\lambda_k^{1})|\cdot
        |g_1(\lambda)|,
\eeqa
so that
\bea\label{estimate4}
  \log\frac{1}{|B_{\Lambda\setminus \{\lambda_j^{1}\}_j}(\lambda_k^{1})|}
        \le P[w_1](\lambda_k^{1}),\quad  k\in\N.
\eea
%Without loss of generality
Replacing possibly $w_1$ by $(w_1)_+$ we can assume $w_1\ge 0$, so that
$P[w_1]$ is a positive harmonic function.
%Since $h_2 T_2$ does not vanish and is bounded above by 1, the function
%$\log |h_2 T_2|$ is a negative harmonic function. 
By Harnack's
inequality, there exists a constant $c_H\ge 1$ such that
\beqa
       \frac{1}{c_H} P[w_1](\lambda_k^{1})
       \le P[w_1](z) \le c_H P[w_1](\lambda_k^{1}),
       \quad z\in Q_k,
\eeqa
%hence
%\beqa
% P[w_1](\lambda_k^{1})|^{c}\le |g_1(z)|
%       \le |g_1(\lambda_k^{1})|^{1/c},\quad z\in Q_k.
%\eeqa
%This yields
%\bea\label{estimate4}
%       |(h_2 T_2)^c(\lambda')|\le |(h_2 T_2)(\lambda_k^{1})|\le
%       |B_{\Lambda\setminus \{\lambda_j^{1}\}_j}(\lambda_k^{1})|
%\eea
and in particular for every $\lambda'\in \Lambda_1\cap Q_k$.

At this point we have to change the argument from \cite{HMNT}.
Instead of using Proposition 4.1 of that paper we have to invoke
Corollary \ref{intsep}.
By construction,
the sequence $\{\lambda_j^{1}\}_j\subset\Lambda_1$ is
separated. Moreover, as a subsequence of an interpolating sequence
it is also interpolating (cf.\ Remark \ref{l-infinit}).
% {\sc or
%the corresponding part in the definition of free interpolation}).
By Corollary \ref{intsep} 
%(see also the reformulation thereafter),
there exists %an outer function $h_1\in \hf$ 
a function $v_1\in L^{\vp}$ such that
\beqa
  \log\frac{1}{|B_{\{\lambda_j^{1}\}_j\setminus \{\lambda_k^{1}\}}
  (\lambda_k^{1})|}
       \le P[v_1](\lambda_k^{1}),\quad k\in\N.
\eeqa
Recall that the weight $v_1$ %defining the outer function $h_1$ 
can be supposed positive, so that $P[v_1]$ is a positive harmonic
function,
%$\log |h_1|$ is a negative harmonic function.
and again by Harnack's inequality
we get
\beqa
       P[v_1](\lambda_k^{1}) \le c_H P[v_1](\lambda')
\eeqa
for every $\lambda'\in \Lambda_1\cup Q_k$.
This together with \eqref{estimate4} and our definition
of $\lambda_k^{1}$ give
\beqa
  \log\frac{1}{|B_{\Lambda\setminus\{\lambda'\}}(\lambda')|}
   &\le& \log\frac{1}{|B_{\Lambda\setminus\{\lambda_k^{1}\}}(\lambda_k^{1})|}
    =\log\frac{1}{|B_{\Lambda\setminus \{\lambda_j^{1}\}_j}(\lambda_k^{1})|}
       +\log\frac{1}{|B_{\{\lambda_j^{1}\}_j\setminus \{\lambda_k^{1}\}}
         (\lambda_k^{1})|}\\
       &\le& P[w_1+v_1](\lambda_k^1)\le P[c_H(w_1+v_1)](\lambda')
\eeqa
for every $\lambda'\in Q_k$ and $Q_k\in Q^{(1)}$. 
We set $u_1:=c_H(w_1+v_1)$ which is clearly in $L^{\vp}$.

Construct in a similar way functions $u_i$
for the sequences $\Lambda_i$, $i=2,3,4$, and set $u=\sum_{i=1}^4u_i$
which is still in $L^{\vp}$ by the quasi-triangular inequality
\eqref{deltaineq}.
So, whenever $\lambda\in\Lambda$, there exists
$k\in \{1,2,3,4\}$ such that $\lambda\in\Lambda_k$, and hence
\begin{equation}
\label{blaschkeminor}
  \log\frac{1}{|B_{\lambda}(\lambda)|}\le 
 P[u_k](\lambda)\le P[u](\lambda).
\end{equation}
%\end{proof}

\section{The trace spaces}\label{trace}

The prove of the trace space characterization is
easier than that given in (see \cite{HMNT}) for the Nevanlinna
class.

In order to see that (c) in Theorem \ref{CNS} implies free interpolation
it suffices to
observe that  $\ell^{\infty}\subset l_{\Phi}$ and use Remark
\ref{l-infinit}.

%Consider the converse.
Assume now that $\Lambda$ is of free interpolation.
%Assume that
Suppose that $(a_\lambda)_\lambda\in \hf|\Lambda$ and $f\in \hf$ is such
that $f(\lambda)=a_\lambda$, $\lambda\in\Lambda$. 
Let $w\in L^1$ be the representing measure of the outer part of 
$f$, i.e.\ $f=If_w$,
%\beqa
% f(z)=I(z)\exp\left(\int_{\T}\frac{\zeta+z}{\zeta-z}w(\zeta)dm(\zeta)
% \right),\quad z\in \DD,
%\eeqa
where $I$ is inner and $w_+\in L^{\vp}$. 
Obvioulsy $\log^+|a_{\lambda}|\le
\log^+(\exp(P[w](\lambda))) 
\le P[w_+](\lambda)$, and so
$(a_{\lambda})\in l_{\Phi}$.
%. Since $w$ can be supposed positive,
%we have $\log^+|a_{\lambda}|

Conversely, suppose that $(a_{\lambda})_{\lambda}$ is
such that there is a positive function $w\in L^{\vp}$ with 
$\log^+|a_{\lambda}|\le P[w](\lambda)$. 
Since $f_w\in \hf$ and $\log^+|f_w|=\log |f_w|=P[w]$ we have
$|a_{\lambda}|\le |f_w(\lambda)|$ for every $\lambda\in \Lambda$.
Since $\Lambda$ is of free interpolation, i.e.\ $\hf|\Lambda$ is
ideal, there exists a function $f_0\in \hf$ interpolating
$(a_{\lambda})_{\lambda}$.
\hfill \qedsymbol

\section{Harmonic majorants}\label{harmmaj}

We begin by recalling the definition of the Poisson balayage: for a
positive finite measure $\mu$ on the closed unit disk we set
\beqa
 B(\mu)(\zeta)=\int_{\DD}\frac{1-|z|^2}{|\zeta-z|^2}\;d\mu(z)
 =\int_{\DD}P_z(\zeta)\,d\mu(z).
\eeqa 
Let $\Har_{\vp}^+=\{P[w]:0\le w \in L^{\vp}\}$. We will begin with 
%the following result that is the
an analog of \cite[Proposition 6.1]{HMNT}
and that does not require the $\nabla_2$-condition.

\begin{proposition}
Let $\vp$ be a strongly convex function, and
let $\mu$ be a positive finite Borel measure on $\T$. Then
$\|B(\mu)\|_{(L^{\vp})^*}<\infty$ if and only if for every $f\in
\Har_{\vp}^+$ we have $\int_{\DD}h\; d\mu<\infty$. Moreover we
have the following relation:
\beqa
 \|B(\mu)\|_{(L^{\vp})^*}=\sup\left\{\int_{\DD}h\; d\mu:
 h=P[w]\in \Har_{\vp}^+, \|w\|_{\vp}\le 1\right\}.
\eeqa
\end{proposition}

So, this proposition furnishes a description of those positive
finite measures on $\DD$ that act against positive harmonic
functions $P[w]$ with $w\in L^{\vp}$.

The proof of this result is short. It is essentially based on an application
of Fubini's theorem and the definition of the norm in 
%the space of 
$L^{\vp}$ by duality. We give it for completeness.

\begin{proof}
%By Fubini's theorem,
\beqa
\lefteqn{ \sup\left\{\int_{\DD}h\; d\mu:
   h=P[w]\in \Har_{\vp}^+, \|w\|_{\vp}\le 1\right\}}\\
 &&=\sup\{\int_{\DD} \int_{\T} \frac{1-|z|^2}{|\zeta-z|^2}
   w(\zeta)\;dm(\zeta)\; d\mu(z):0\le w\in L^{\vp}, \|w\|_{\vp}\le 1\}\\
 &&=\sup\{\int_{\T} w(\zeta) \int_{\DD} \frac{1-|z|^2}{|\zeta-z|^2}\;d\mu(z)
  \;dm(\zeta) :0\le w\in L^{\vp}, \|w\|_{\vp}\le 1\}\\
 &&=\sup\{\int_{\T} w(\zeta) B(\mu)(\zeta)
  \;dm(\zeta) : w\in L^{\vp},\|w\|_{\vp}\le 1\}\\
 &&=\|B(\mu)\|_{(L^{\vp})^*}=\|B(\mu)\|_{L^{\vp^*}}.
\eeqa
%In the last equality w
In the above identities,
we have also used the fact that $\mu$ is a positive
measure so that its balayage is also positive. Hence it is enough to
test against positive functions in $L^{\vp}$.
\end{proof}

We are now in a position to prove Theorem \ref{charharmmaj}.

\begin{proof}[Proof of Theorem \ref{charharmmaj}]
Condition (b) is clearly necessary. It suffices indeed to plug the
estimate $u\le P[w]$ into the chain of equalities in the previous
proof: if $u\le h=P[w]$ for some $0\le w\in L^{\vp}$
then $0\le \int_{\DD}u d\mu\le \int_{\DD}hd\mu\le \|w\|_{\vp}
\|B\mu\|_{\vp^*}$ (H\"older's inequality for Orlicz spaces,
see \cite[Theorem 9.3]{KR}).

Let us consider the sufficiency. So suppose that
there exists a constant $C\ge 0$ such that
\bea\label{condmaj}
 \sup_{\mu\in {\mathcal B}_{\vp^*}} \int u(z)\,d\mu(z)\le C,
\eea
where ${\mathcal B}_{\vp^*}=\{\mu :$ positive measure on $\DD$ such that
$\|B\mu\|_{(L^{\vp})^*}\le 1\}$.
We want to prove that $u$ admits a harmonic majorant $P[w]$ with
$w\in L^{\vp}$.

We will begin as in \cite{HMNT} 
by discretizing the problem (see in particular \cite[Lemma 6.3]{HMNT}). 
For this, let again $Q_{n,k}$ be the
dyadic cubes and $z_{n,k}$ the corresponding center.
Fix $z_{n,k}^*\in Q_{n,k}$ such that $u(z_{n,k}^*)\ge (\sup_{Q_{n,k}}u)/2$.
%Define
%\beqa
% \hat{u}(z)=
% \left\{\begin{array}{l}
%  \sup_{z\in Q_{n,k}} u(z) \text{ if }z=z_{n,k}^*, \\
%  0 \text{ otherwise}.
% \end{array}
% \right.
%\eeqa
We will set $u_{n,k}:={u}(z_{n,k}^*)$ and $\hat{u}_{n,k}=\sup_{Q_{n,k}}u$.
Let us check that if $u$ satisfies \eqref{condmaj} then
there is a constant $C'$ such that whenever
$(c_{n,k})$ is a finite sequence of non-negative coefficients with
\beqa
 \|\sum c_{n,k}P_{z_{n,k}}\|_{\vp^*}\le 1,
\eeqa
we get
\beqa
 \sum c_{n,k} %\sup_{Q_{n,k}}{u}
 \hat{u}_{n,k}\le C'.
\eeqa

Indeed, setting $\mu:=\sum c_{n,k}\delta_{z_{n,k}^*}$ we obtain a
positive finite measure on $\DD$ such that
\beqa
 \|B(\mu)(\zeta)\|_{\vp^*}=\|\int_{\DD}P_z(\zeta)\,d\mu(z)\|_{\vp^*}
 = \|\sum c_{n,k}P_{z_{n,k}^*}\|_{\vp^*}
 \le K \|\sum c_{n,k}P_{z_{n,k}}\|_{\vp^*}
 \le K.
\eeqa
In the last estimate we have used the existence of  a
constant $K=K(\delta)$ such that  $|b_u(v)|\le \delta$
implies that $\frac{1}{K}P_v(\zeta)\le  P_u(\zeta)\le K P_v(\zeta)$ 
for all $\zeta\in \T$
and the fact that the Orlicz space $L^{\vp}$ has the lattice
property (which can be seen by using the dual representation
of the norm).

Now, using \eqref{condmaj}, we get
\beqa
 0\le\sum c_{n,k} \hat{u}_{n,k} %= \sum c_{n,k}\sup_{Q_{n,k}}u
 \le 2\sum c_{n,k}u_{n,k}=2\int {u}d\mu \le 2CK
\eeqa
So
\beqa
 \sup\left\{\sum c_{n,k} \hat{u}_{n,k}:c_{n,k}\ge 0 
 \text{ for all }n,k\text{ and }
 \|\sum c_{n,k}P_{z_{n,k}}\|_{\vp^*}\le 1\right\}\le 2CK,
\eeqa
in other words, for all positive finite sequences $(c_{n,k})$ we have
\beqa
 \sum c_{n,k} \hat{u}_{n,k}\le 2CK \|\sum c_{n,k}P_{z_{n,k}}\|_{\vp^*}.
\eeqa
%Let $E_0=\{\sum c_{n,k} P_{z_{n,k}^*}:(c_{n,k})$ finite sequence$\}$ 
%equipped with the
%norm $\|\cdot\|_{\vp^*}$ and
Let $E$ be the closure in $L^{\vp^*}$ of the
$\R$-space generated by the Poisson kernels
$P_{z_{n,k}}$ for all $n,k$ (which in fact corresponds to $L^{\vp^*}$).
By the theorem of Mazur-Orlicz 
(see e.g.\ \cite[Chapter 2, Proposition 2.2]{Pe}),
there exists a linear continuous mapping $T:E\lra \R$ with same norm $2CK$
such that 
\beqa
 \hat{u}_{n,k}\le T P_{z_{n,k}}.
\eeqa
So $T\in E^*=(L^{\vp^*})^*=L^{\vp^{**}}=L^{\vp}$. Hence there
exists $w_u\in L^{\vp}$ such that for all $v=\sum c_{n,k}
P_{z_{n,k}}$ we have
$Tv=\int_{\T} vw_u dm=\sum c_{n,k}P[w_u](z_{n,k})$. In particular
for every $n,k$
\beqa
 0\le \hat{u}_{n,k}\le TP_{z_{n,k}}= P[w_u](z_{n,k})\le P[(w_u)_+](z_{n,k}).
\eeqa
By Harnack's inequality the last term is bounded by 
$P[c_H(w_u)_+](z)$ for % some constant $c\ge 1$
and every $z\in Q_{n,k}$ so that for every $z\in
Q_{n,k}$
\beqa
 0\le u(z)\le \hat{u}_{n,k}\le P[(w_u)_+](z_{n,k})\le P[c_H(w_u)_+](z).
\eeqa
Since this is true for all $n,k$ and since 
clearly $c_H(w_u)_+\in L^{\vp}$ we have achieved the proof.
\end{proof}
%\qed

%Let $E_0=\{\sum c_{n,k} P_{z_{n,k}^*}:(c_{n,k})$ finite sequence$\}$ 
%equipped with the
%norm $\|\cdot\|_{\vp^*}$ and $E$ its closure in $L^{\vp^*}$
%(which in fact corresponds to $L^{\vp^*}$). Let us consider the
%linear (positive) application
%\beqa
% \psi_u : E_0&\lra& \R,\\
%     \sum c_{n,k} P_{z_{n,k}^*} &\lmto& \sum c_{n,k} \hat{u}_{n,k}.
%\eeqa

%$\hat{u}$ (with possibly different constants). It is clearly
%enough to check the condition for measures supported on 
%$(z_{n,k})_{n,k}$. Also we can restrict ourselves to finitely
%supported measures

\section{An example}\label{examples}

In this section we will consider concrete separated sequences and
check whether they are interpolating for Hardy-Orlicz spaces $\hf$
associated with $\Phi=\Phi_{\eps}=\psi_{\eps}\circ\log$ where
\beqa
 \psi_{\eps}=t\log^{\eps}t
\eeqa
for some $\eps>0$ (even if this has no special meaning for our
situation one could observe that for $\eps=1$ the space $L^{\psi_1}$ 
is the Zygmund space $L\log L$).
%%Our first example is 
%We construct a separated sequence that is not
%interpolating for any $\hf$ where $\vp(t)\ge t\log^{\eps} t$ for 
%$t\ge t_0$ and $\eps>0$. %, for some $p>1$.

\begin{proposition}
For every $\eps>0$ there exists a separated sequence $\Lambda$ that
is interpolating for ${\mathcal H}_{\psi_{\delta}}$ whenever 
$0<\delta<\eps$ but not for ${\mathcal H}_{\psi_{\eps}}$.
\end{proposition}

\begin{proof}
Fix $\eps>0$ and
let $\lambda_{n,k}=(1-1/2^n)e^{k2\pi i/2^n}$, where $n\in\N^*$,
$k\in \{-k_n,\ldots,k_n\}$ and $k_n:=[2^n/(n\log^{1+\eps} n)]$. Then
$\sum_{n,k}(1-|\lambda_{n,k}|)\simeq\sum_n [1/(n\log^{1+\eps}n)]<\infty$
so that $\Lambda=\{\lambda_{n,k}\}_{n,k}$ satisfies the Blaschke
condition.

Since $\Lambda$ is separated we have
for $\lambda,\mu\in\Lambda$ by standard estimates
\beqa
 \log\frac{1}{|b_{\mu}(\lambda)|} \simeq 1-|b_{\mu}(\lambda)|
 \simeq \frac{(1-|\mu|)(1-|\lambda|)}{|1-\overline{\mu}\lambda|^2}.
\eeqa
We will compute $\log(1/|B_{\lambda}(\lambda)|)$ for
$\lambda=\lambda_{n,0}=1-1/2^n$, $n\in\N$, and show that this
exceeds the maximal admissible growth in ${\mathcal H}_{\Phi_{\eps}}$.
%Set $k_n:=[2^n/(n\log^2n)$.
%\load{\scriptsize}{\scs}
\beqa
 \log\frac{1}{|B_{\lambda}(\lambda)|}
 &=&\sum_{j\ge
   1}\sum_{
 \begin{array}{c}
 {\scriptstyle l=-k_j}\\
 {\scriptstyle l\neq 0\text{ if }j=n}
 \end{array}}^{k_j}\log\frac{1}{|b_{\lambda_{j,l}}(\lambda)|}
 \simeq \sum_{j\ge 1}\sum_{
 \begin{array}{c}
 {\scriptstyle l=-k_j}\\
 {\scriptstyle l\neq 0\text{ if }j=n}
 \end{array}}^{k_j}
   \frac{(1-|\lambda|)(1-|\lambda_{j,l}|)}{|1-\overline{\lambda_{j,l}}
   \lambda|^2}\\
 &=&\frac{1}{2^n}\sum_{j\ge 1}\frac{1}{2^j}\sum_{
 \begin{array}{c}
 {\scriptstyle l=-k_j}\\
 {\scriptstyle l\neq 0\text{ if }j=n}
 \end{array}}^{k_j}
  \frac{1}{|1-\overline{\lambda_{j,l}}
   \lambda|^2}\\
 &\ge& \frac{1}{2^n}\sum_{j\ge 2n}\frac{1}{2^j}\sum_{l=-k_j}^{k_j}
  \frac{1}{|1-\overline{\lambda_{j,l}}
   \lambda|^2}
\eeqa
%Note that
We get for $j\ge 2n$
\bea\label{eq5}
 \frac{1}{2^j}\sum_{l=-k_j}^{k_j}
  \frac{1}{|1-\overline{\lambda_{j,l}}
   \lambda|^2}=
 \sum_{l=-k_j}^{k_j} \frac{1}{2^j}\frac{1}{|e^{l2\pi i/2^j}-r_{n,j}|^2},
\eea
where $r_{n,j}=(1-1/2^n)(1-1/2^j)$. It can be noted that for fixed
$n$, 
$r_{n,j}$ goes rapidly and increasingly to $\lambda=1-1/2^n$
as $j\to+\infty$
(see Figure 2). The sum in \ref{eq5} is a Riemann
sum for 
\beqa
 \int_{I_j} \frac{1}{|e^{i\vt}-r_{n,j}|^2}\; d\vt,
\eeqa
where $I_j=[-2\pi k_j/2^j,2\pi k_j/2^j]=
[-2\pi/(j\log^{1+\eps} j),2\pi/(j\log^{1+\eps} j)]$
(note that the intervalle $I_j$ on which we integrate depends
on $j$). 
Then, since the function to be
integrated is continuous, we get for fixed $n$
and $j\ge 2n$
\beqa
 \frac{1}{2^j}\sum_{l=-k_j}^{k_j}
  \frac{1}{|1-\overline{\lambda_{j,l}}
   \lambda|^2}\simeq \int_{I_j} \frac{1}{|e^{i\vt}-r_{n,j}|^2}\; d\vt
%,\quad j\to +\infty
\eeqa

%(note that we are interested in $j\ge 2n$).
(for a better control on the constants it is possible to replace
$j\ge 2n$ by $j\ge 2n+K$ for some fixed $K>0$).

\begin{center}
\begin{picture}(0,0)%
\includegraphics{noypoiss.pstex}%
\end{picture}%
\setlength{\unitlength}{2368sp}%
\begingroup\makeatletter\ifx\SetFigFont\undefined%
\gdef\SetFigFont#1#2#3#4#5{%
  \reset@font\fontsize{#1}{#2pt}%
  \fontfamily{#3}\fontseries{#4}\fontshape{#5}%
  \selectfont}%
\fi\endgroup%
\begin{picture}(4524,3843)(2689,-6892)
\put(5026,-3286){\makebox(0,0)[lb]{\smash{{\SetFigFont{7}{8.4}{\familydefault}{\mddefault}{\updefault}{\color[rgb]{0,0,0}$r_{n,j}$}%
}}}}
\put(3901,-6361){\makebox(0,0)[lb]{\smash{{\SetFigFont{7}{8.4}{\familydefault}{\mddefault}{\updefault}{\color[rgb]{0,0,0}$\underbrace{\phantom{andreashartmann}}_{ }$}%
}}}}
\put(4651,-6811){\makebox(0,0)[lb]{\smash{{\SetFigFont{10}{12.0}{\familydefault}{\mddefault}{\updefault}{\color[rgb]{0,0,0}$I_j$}%
}}}}
\put(4951,-3886){\makebox(0,0)[lb]{\smash{{\SetFigFont{7}{8.4}{\familydefault}{\mddefault}{\updefault}{\color[rgb]{0,0,0}$\lambda=1-1/2^n$}%
}}}}
\put(4876,-4411){\makebox(0,0)[lb]{\smash{{\SetFigFont{7}{8.4}{\familydefault}{\mddefault}{\updefault}{\color[rgb]{0,0,0}$\alpha_{n,j}$}%
}}}}
\put(2776,-6361){\makebox(0,0)[lb]{\smash{{\SetFigFont{10}{12.0}{\familydefault}{\mddefault}{\updefault}{\color[rgb]{0,0,0}$\T$}%
}}}}
\end{picture}%
\label{figure3}\\
%\vspace{0.5cm}
Figure 2 
\end{center}

It is well known (see e.g.\ \cite[p.13]{Gar}) that 
\beqa
 \int_{I_j} \frac{1-r_{n,j}^2}{|e^{i\vt}-r_{n,j}|^2}\; d\vt
 \simeq \frac{\alpha_{n,j}}{\pi}
\eeqa
($\alpha_{n,j}$ is the angle indicated in Figure 2).
As we have already mentioned $1-r_{n,j}\ge 1-\lambda=1/2^n$,
and so for fixed $n$ we get
$\alpha_{n,j}\to 0$ as $j\to+\infty$. 
So $\alpha_{n,j}=2\alpha_{n,j}/2\sim 2 \tan (\alpha_{n,j}/2)=
2 (|I_j|/2)/(1-r_{n,j})
\sim|I_j|/(1-\lambda)=2^n |I_j|=4\pi 2^n/(j\log^{1+\eps}j)$.
Hence
\beqa
 \log\frac{1}{|B_{\lambda}(\lambda)|}
 &\ge&\frac{1}{2^n}\sum_{j\ge 2n}\frac{1}{1-r_{n,j}^2}
  \sum_{l=-k_j}^{k_j} \frac{1}{2^j}
   \frac{1-r_{n,j}^2}{|e^{l2\pi i/2^j}-r_{n,j}|^2}\\
 &\simeq& \frac{1}{2^n}\sum_{j\ge 2n}\frac{1}{1-r_{n,j}^2}
  \int_{I_j} \frac{1-r_{n,j}^2}{|e^{i\vt}-r_{n,j}|^2}\; d\vt\\
 &\simeq&\frac{1}{2^n}\sum_{j\ge 2n}\frac{1}{1-r_{n,j}^2}
   \frac{4 \cdot 2^n}{j\log^{1+\eps}j}\\
 &=& \sum_{j\ge 2n}\frac{1}{1-r_{n,j}^2}\frac{4}{j\log^{1+\eps}j}
\eeqa
Moreover $j\ge 2n$, so that $1-r_{n,j}^2\simeq 1-r_{n,j}
=\frac{1}{2^n}+\frac{1}{2^j}-
\frac{1}{2^{n+j}}\simeq \frac{1}{2^n}$. Using  $1-|\lambda|=1-\lambda
=1/2^n$ and so $n= \log(1/(1-|\lambda|))/\log 2$ we get
\bea\label{estimexemple}
 \log\frac{1}{|B_{\lambda}(\lambda)|}
 &\gtrsim&  2^n  \sum_{j\ge 2n}\frac{4}{j\log^{1+\eps}j}
 \sim  2^{n} \frac{4}{\log^{\eps} (2n)}\nn\\
 &\sim& \frac{4}{1-|\lambda|}\frac{1}{\log^{\eps} n}\sim
 \frac{4}{1-|\lambda|}\frac{1}{\log^{\eps}\log \frac{1}{1-|\lambda|}}.
\eea

Let us check that this is not compatible with
free interpolation in ${\mathcal H}_{\Phi_{\eps}}$.
Suppose to the contrary that $\Lambda$ is an interpolating sequence
for ${\mathcal H}_{\Phi_{\eps}}$. Then,
by Theorem \ref{CNS},
there exists a positive function $w\in L^{\psi_{\eps}}(\T)$ such that 
$\log(1/|B_{\lambda}(\lambda)|)\le P[w](\lambda)$ for every
$\lambda\in\Lambda$. 
As we have already discussed in 
Section \ref{HO},  $f(z)=\exp(\int (\zeta +z)/(\zeta-z) w(\zeta)
\;dm(\zeta))$ is an outer function in ${\mathcal H}_{\Phi_{\eps}}$. Hence
\beqa
 P[w](z)=\log|f(z)|\le \psi_{\eps}^{-1}\left(\frac{c_f}{1-|z|}\right).
\eeqa
Note that $\psi_{\eps}(u/\log^{\eps} u)\sim u$ as $u\to +\infty$,
so that
$\psi_{\eps}^{-1}(u)\sim u/\log^{\eps} u$. Hence, setting $u=1/(1-|z|)$
we get for $z$ sufficiently close to $\T$,
\beqa
  P[w](z) \lesssim c_f \frac{1}{1-|z|}\frac{1}{\log^{\eps} \frac{c_f}{1-|z|}}.
\eeqa
Since the right hand side of the last inequality is negligible
with respect to the right hand side of \eqref{estimexemple}
we have reached a contradiction.
%$\log(1/|B_{\lambda}(\lambda)|)$ has no quasi-bounded harmonic
%majorant with weight in $L^{\psi_{\eps}}$. So, by Theorem \ref{CNS}
So, the sequence $\Lambda$ is not interpolating for ${\mathcal H}_{\Phi_{\eps}}$
and hence for no ${\mathcal H}_{\Phi}$ with
$\Phi=\vp\circ\log$ and $\vp$ a stronly convex function with
$\vp(t)\ge \psi_{\eps}(t)$, $t\ge t_{\eps}$.

In order to finish the proof, we check that
the above constructed sequence is interpolating for 
${\mathcal H}_{\Phi_{\delta}}$ whenever $0<\delta<\eps$.
Since $\Lambda$ is separated it is of course interpolating for the
Smirnov (and Nevanlinna) class, see \cite[Corollary 1.9]{HMNT},
which means that there is a function $u\in L^1(\T)$ such that
\bea\label{majQB}
 \log(1/|B_{\lambda}(\lambda)|)\le P[u](\lambda).
\eea
It is known that the function $u$ can be chosen explicitely by:
\beqa
 u=c_0\sum_{\lambda\in\Lambda}\chi_{I_{\lambda}}
\eeqa
(see \cite[Proposition 4.1]{HMNT} and also
\cite[p.124]{NPT}),
the intervalle $I_{\lambda}=\{e^{it}\in\T:
|t-\arg \lambda|\le c(1-|\lambda|)\}$ appearing in the above
formula being the so-called Privalov shadow.
It turns out that in the present situation this intuitive candidate
for $u$ is the right one to get a harmonic majorant. In other words,
we have to check that $u\in L^{\psi_{\delta}}$, $0<\delta<\eps$, 
and this will finish the proof. 
%This is what we will heck now.
So, let us suppose that the constant
$c$ in the definition of $I_{\lambda}$ is adapted in such a way that
$I_{\lambda_{n,k}}$ and $I_{\lambda_{n,k+1}}$ touch without overlap
(this is not really of importance). We then consider the shadow of
the stage $n$: $\bigcup_{j=1,...,k_n}I_{\lambda_{n,j}}=
[-1/2^n-1/(n\log^{1+\eps}n),1/(n\log^{1+\eps}n)+1/2^n]$ which is essentially
the interval $[-1/(n\log^{1+\eps}n),1/(n\log^{1+\eps}n)]$. So the function $u$
is essentially equal to $k$ on $ [-1/(k\log^{1+\eps}k),
-1/((k+1)\log^{1+\eps}(k+1))[ \cup ]1/((k+1)\log^{1+\eps}(k+1)),
1/(k\log^{1+\eps}k)]$. In order that $\psi_{\delta}\circ u\in L^1$ it
thus suffices (using some Fubini) that
\beqa
 \sum_{k\ge 1} \frac{\psi_{\delta}(k+1)-\psi_{\delta}(k)}
 {k\log^{1+\eps}k}<\infty,
\eeqa
and this holds for $0<\delta<\eps$.
\end{proof}

If one wishes to get closer to $L^1$, one could e.g.\ consider
strongly convex functions $\vp(t)=t\log^{\eps}\log t$ by
choosing $k_j=2^j/(j\log j\log^{1+\eps}\log j)$.

We wanted to emphasize in this section on the behaviour of 
{\it separated} Blaschke
sequences since these were already interpolating for the Smirnov
class. Our examples make clear that the situation is much more
delicate in big Hardy-Orlicz spaces. 

Another and of course easier
way of producing examples of interpolating sequences for
our spaces is to take two $\Hi$-interpolating sequences that approach
in a critical way: if $\Lambda_1=\{\lambda_n\}_n$ 
is such an $\Hi$-interpolating sequence, take $\mu_n\neq\lambda_n$ close to
$\lambda_n$ and define $\Lambda_2=\{\mu_n\}_n$. 
Then, for $B=B_{\Lambda_1\cup\Lambda_2}$ we get 
$\log (1/|B_{\lambda_n}(\lambda_n)|)\simeq \log
(1/|b_{\lambda_n}(\mu_n)|=:\eta_n$, so that suitable choices of
$(\eta_n)_n$ yield interpolating sequences for
some big Hardy-Orlicz spaces which are not interpolating for others
(see e.g.\ \cite[Example 9.2]{HMNT} and \cite{HM2} for such
constructions).

{\bf Acknowledgements:} K. Dyakonov
reminded me that Theorem \ref{hoffthm} is due to Hoffman.

\providecommand{\bysame}{\leavevmode\hbox to3em{\hrulefill}\thinspace}
\providecommand{\MR}{\relax\ifhmode\unskip\space\fi MR }

\end{document}